\documentclass[12pt]{article}
\title{Symplectic Resolutions for Symmetric Products of Surfaces}
\author{Baohua Fu}
\usepackage{amsmath,amssymb,amsthm, amscd}
\chardef\bslash=`\\
\newtheorem{Thm}{Theorem}

\newtheorem{Lem}[Thm]{Lemma}

\newtheorem{Rque}{Remark}

\def\cit{{\mathbb C}}

\def\qit{{\mathbb Q}}

\def\0{{\mathcal O}}

\def\m{{\mathfrak m}}

\begin{document}
\maketitle
\begin{abstract}
Let $S$ be a smooth complex connected analytic surface which admits a holomorphic symplectic structure. Let $S^{(n)}$ be its $n$th symmetric product.
We prove that every projective symplectic resolution of $S^{(n)}$ is isomorphic to the Douady-Barlet resolution $S^{[n]} \rightarrow S^{(n)}$.
\end{abstract}

\section{Introduction}

Let $X$ be a normal complex analytic variety. Following \cite{Bea}, the variety $X$  
is said to have {\em symplectic singularities} if there 
exists a holomorphic
symplectic 2-form $\omega$ on $X_{reg}$ such that for any resolution of singularities $\pi: \widetilde{X} \rightarrow X$, the 2-form
$\pi^* \omega$ defined a priori on $\pi^{-1}(X_{reg})$ can be extended to a holomorphic 2-form on $\widetilde{X}$. If furthermore 
the 2-form $\pi^* \omega$ extends to a holomorphic symplectic 2-form on the whole of $\widetilde{X}$ for some resolution of $X$, 
then we say that $X$ admits a {\em symplectic resolution}, and the resolution $\pi$ is called {\em symplectic}.

There are two classes of examples of symplectic singularities. One consists of  normalizations of  closures of  nilpotent
orbits in a complex semi-simple Lie algebra. For these singularities, we proved in \cite{Fu} that every symplectic resolution is isomorphic
to the collapsing of the zero section of the cotangent bundle of a projective homogeneous space.

The other class of examples  consists of so called quotient singularities, i.e. 
singularities of the form $V/G$, with $V$ a symplectic variety and $G$ a finite group of symplectic 
automorphisms of $V$. 
 In this note, we are interested in the following special case. Let $S$ be a complex analytic manifolds which admits a holomorphic symplectic
structure. Then by Proposition 2.4 (\cite{Bea}), the symmetric product $S^{(n)} = S^n/\Sigma_n$ is a normal variety with symplectic
singularities, where $\Sigma_n$ is the permutation group of $n$ letters. By Theorem 2.1 (\cite{Fub}), $S^{(n)} $ (for $n > 1$)  admits no 
symplectic resolution
as soon as $dim(S) > 2$. From now on, we will suppose that $dim(S) = 2$, i.e. $S$ is a smooth complex connected analytic surface admitting a symplectic 
structure. Examples of such surfaces include $K3$ surfaces, abelian surfaces, cotangent bundles of algebraic curves etc..

 Recall that the Douady space $S^{[n]}$  parametrizes zero-dimensional analytic
 subspaces of length $n$ in $S$. It is well-known that $S^{[n]}$ is a $2n$-dimensional smooth complex manifolds. If $S$ is algebraic, it is the usual
Hilbert scheme $Hilb^n(S)$ of points on $S$. The  Douady-Barlet morphism $S^{[n]} \rightarrow S^{(n)}$ provides a projective resolution of singularities
for $S^{(n)}$. If $S$ admits a symplectic structure, then $S^{[n]}$ is again symplectic (see \cite{Be}). In this case 
$S^{[n]}$ gives a projective symplectic resolution for $S^{(n)}$. Our purpose of this note is to prove the following:
\begin{Thm}
Let $S$ be a smooth connected complex  analytic surface, which admits a symplectic structure. Then every projective symplectic
resolution of $S^{(n)}$ is isomorphic to the Douady-Barlet resolution $S^{[n]} \rightarrow S^{(n)}$.
\end{Thm}

One should bear in mind that in general there may be more than one symplectic resolutions for a variety with symplectic singularities.
Such an  example is given by a Richardson nilpotent orbit admitting two or more polarizations with non-conjugate Levi factors (see \cite{Fu}).

A similar result for symplectic quotients of $\cit^{2n}$ has been proved by D. Kaledin (see Theorem 1.9 \cite{Ka}).

{\em Acknowledgments:} I would like to thank A. Beauville and A. Hirschowitz for helpful discussions.

\section{Proof of the theorem}
\begin{Lem}\label{lem1}
Let $X$ be a normal variety with symplectic singularities and $U$ its smooth part. Let $Sing(X)=\cup_{i=1}^k F_i$ be 
the decomposition into irreducible components of the singular part of $X$. Suppose that $X$ admits a projective 
symplectic resolution  $\pi:Z \rightarrow X$.
 Furthermore we suppose that: 

(i). $Pic(U)$ is a torsion group;

(ii). for any $i$, there exists an analytic proper sub-variety $B_i$ in $F_i$, such that the restricted map $\pi': Z_* \rightarrow X_*$
is the blow-up of $X_*$ with center $\cup_i (F_i-B_i)$, where $X_* = X - \cup_i B_i$ and $Z_* = \pi^{-1}(X_*)$.

Then  every $F_i$ is of codimension 2 in $X$, and the morphism $\pi: Z \rightarrow X$ is isomorphic to the blow-up 
$Bl(X, \cap_i \m_i^{d_i})$ for some integers $d_i$, where $\m_i$ is the  ideal $\m_{F_i}$ defining $F_i$.
\end{Lem}
\proof
That $X$ is normal implies $X-U$ has codimension at least $2$ in $X$, thus the Weil divisor group $Cl(X)$ is isomorphic to 
the Picard group $Pic(U)$, which is of torsion by hypothesis (i). This shows that $X$ is $\qit$-factorial. Now the first affirmation follows 
from Corollary 1.3 \cite{Fu}.

We will use an idea of D. Kaledin \cite{Ka} to prove the second affirmation. Since $\pi$ is projective, $Z = Proj \oplus_k \pi_*L^k$
for some holomorphic line bundle $L$ on $Z$.  Notice that $Pic(\pi^{-1}(U)) = Pic(U)$ is of torsion, replacing $L$ by some positive power, 
we can suppose $L|_{\pi^{-1}(U)}$ is trivial, thus $\pi_*L|_U \simeq \0_X|_U.$ Since $X$ is normal and $X-U$ has codimension $\geq 2$,
this gives an embedding $\pi_*L \rightarrow \0_X$, thus $\pi$ is identified with the blow-up $Bl(X, \pi_*L).$

Let $i: Z_* \rightarrow Z$ and $j: X_* \rightarrow X$ be the natural inclusions. We have the following commutative diagram:
\begin{equation*} \begin{CD}
Z_* @>i>> Z \\ @V\pi'VV  @VV{\pi}V \\
X_* @>j>> X.
\end{CD} \end{equation*}
The projection formula gives $L = i_* i^* L$, so
$\pi_* L = \pi_* i_* i^* L = j_* \pi'_* i^*L$. By hypothesis (ii), we have $\pi': Z_* \rightarrow X_*$ is the blow-up of $X_*$ with center
$\cup_i (F_i - B_i)$, thus $\pi'_* i^*L  = j^* (\cap_i \m_i^{d_i})  $ for some positive  integers $d_i$. This gives that 
$\pi_* L =j_* j^* (\cap_i \m_i^{d_i}) = \cap_i \m_i^{d_i}$, which concludes the proof.$\clubsuit$

\begin{Lem}\label{lem2}
Let $\Delta = \{ (z_1, z_2) \in \cit^2 | |z_1| + |z_2| < 1\}$ be the disc in $\cit^2$ and
 $X = \Delta^{(k_1)} \times \cdots \times \Delta^{(k_l)}.$  Then $X$ is a normal variety with symplectic singularities, $Sing(X)$
consists of  $N= \#\{i | k_i > 1\} $ co-dimension 2 irreducible components $F_i$ and any projective symplectic resolution of $X$ is
 isomorphic to a blow-up
$Bl(X, \cap_i \m_i^{d_i})$ for some integers $d_i$, where $\m_i = \m_{F_i}$.
\end{Lem}
\proof
The first affirmation is easy. For the second affirmation, one notices that for $k>1$,  $Sing(\Delta^{(k)})$ is irreducible and of codimension 2 in 
$\Delta^{(k)}$. We will apply Lemma \ref{lem1} to prove the third affirmation. Denote by 
$$p: Y: = \Delta^{k_1} \times \cdots \times \Delta^{k_l}  \rightarrow  X= \Delta^{(k_1)} \times \cdots \times \Delta^{(k_l)}$$
the natural quotient by the product of permutation groups $\Sigma_{k_1} \times \cdots \times \Sigma_{k_l}. $ Let $U$ be the smooth part
of $X$ and $V = p^{-1}(U)$, then $Y - V$ is of codimension 2 in $Y$. In particular, we have $Pic(V) = Pic(Y) = 0$, which gives that 
$Pic(U)$ is of torsion, so  condition (i) of Lemma \ref{lem1} is verified.

We consider one  component, say  $F_1 = Sing(\Delta^{(k_1)}) \times \Delta^{(k_2)} \times \cdots \times \Delta^{(k_l)}$ ($k_1 > 1$).
Denote by $F_{1,*}$ the open set of $F_1$ consists of cycles of the form $2x_1 + x_2 + \cdots + x_{n-1}, $ where $n =\sum_j k_j$,
$x_i$ are distinct points and $2x_1+ x_2+\cdots + x_{k_1-1}$ is in $Sing(\Delta^{(k_1)})$.   We define $B_i = F_i - F_{i,*}$.
In fact, $F_{1,*}$ is the smooth part of $F_1$ and $B_1$ is its singular part.
 Then locally on  $F_{1,*}$, the singularities looks like $\cit^{2n-2} \times (Q,0)$, where 
$Q$ is the cone over a smooth conic and $0$ is its vertex. Recall that every crepant resolution of $Q$ is the blow-up of $Q$ with center  0. 
This gives  that locally on $F_{1,*}$, any symplectic resolution $\pi: Z \rightarrow X$ is the 
blow-up of $X$ along the subvariety $F_{1,*}$, so the second condition of Lemma \ref{lem1} is also verified. 
$\clubsuit$

\proof[End of the proof]
Let $\pi: Z \rightarrow S^{(n)}$ be a projective symplectic resolution. 
For any 0-cycle $I = \sum_i k_ix_i $ in $S^{(n)}$, where $x_i$ are distinct points, we take a small disc $\Delta_i$ near $x_i$ 
on $S$, then $X = \Delta_1^{(k_1)} \times \cdots \times \Delta_l^{(k_l)}$ gives an open neighborhood of $I$ in $S^{(n)}$, so by 
Lemma \ref{lem2}, the restricted map $\pi^{-1}(X) \rightarrow X$ is isomorphic to the blow-up $Bl(X, \cap_i \m_i^{d_i}) $ for some
integers $d_i$. 

Notice that  $F = Sing(S^{(n)})$ is irreducible, so
$F - Sing(F)$ is connected, which implies that all the integers $d_i$ should be equal, since they are locally constant.
  Thus $Z$ is isomorphic to the blow-up 
$Bl(S^{(n)}, \m_F^d)$ for some $d$, which is also isomorphic to  $Bl(S^{(n)}, \m_F)$, i.e. there exists, up to isomorphisms,
 at most one projective
symplectic resolution for $S^{(n)}$, which proves the theorem.
$\clubsuit$

\begin{Rque}
The author does not know whether the product  $S^{(n_1)} \times S^{(n_2)} \times \cdots \times S^{(n_l)}$ admits other projective
symplectic resolutions  than the one given by $S^{[n_1]} \times S^{[n_2]} \times \cdots \times S^{[n_l]}$.
\end{Rque}

\quad \\
Labortoire J.A.Dieudonn\'e, Parc Valrose \\ 06108 Nice cedex 02, FRANCE \\
baohua.fu@polytechnique.org

\begin{thebibliography}{10}
\bibitem[Be]{Be}
A. Beauville, \emph{Vari\'et\'es K\"ahleriennes dont la premi\`ere classe de Chern est nulle}, J. Diff. Geom. {\bf 18}(1983), 755-782
\bibitem[Bea]{Bea}
A. Beauville, \emph{Symplectic Singularities}, Invent. Math. {\bf 139}(2000), 541--549
\bibitem[Fu]{Fu}
B. Fu, \emph{Symplectic resolutions for nilpotent orbits}, Invent. Math. {\bf 151}(2003), 167-186
\bibitem[Fub]{Fub}
B. Fu, \emph{Symplectic resolutions for quotient singularities}, preprint math.AG/0206288
\bibitem[Ka]{Ka}
D. Kaledin, \emph{Dynkin diagrams and crepant resolutions of quotient singularities}, preprint math.AG/9903157, to appear in Selecta Math.

\end{thebibliography}
\end{document}